\documentclass[11pt,titlepage]{article}
\usepackage{amsmath,amsthm, amssymb}
\usepackage{multicol}
\usepackage{color}
\usepackage{secdot}
\allowdisplaybreaks[4]
\makeatletter
\def\section{\@startsection{section}{1}{\z@}{-4ex plus -1ex minus -.2ex}{1.5ex plus .2ex}{\Large\bf}}
\makeatother
\makeatletter

\@addtoreset{equation}{section}     
\makeatother

\def\because{\raisebox{1.1ex}{.}.\raisebox{1.1ex}{.}}

\def\R{\mathbf{R}}

\newtheorem{proposition}{Proposition}
\newtheorem{theorem}{Theorem}

\newtheorem{corollary}{Corollary}
\theoremstyle{definition}

\makeatletter
\def\@proofcounterend{.}
\def\@proofcounterend{}
\def\proof{\@ifnextchar[{\@yproof}{\@xproof}}
\def\@xproof{\@beginproof\ignorespaces}
\def\@yproof[#1]{\@opargbeginproof{#1}\ignorespaces}
\def\@beginproof{\pushQED{\qed}
\topsep6\p@\@plus6\p@\relax
\begin{trivlist}
  \item[\hskip\labelsep{\bf �ؖ�\@proofcounterend}]\rm}
\def\@opargbeginproof#1{\pushQED{\qed}
\topsep6\p@\@plus6\p@\relax
\begin{trivlist}
  \item[\hskip\labelsep{\bf #1�̏ؖ�\@proofcounterend}]\rm}
\def\endproof{\popQED \end{trivlist}}  
\def\subproof{\@beginsubproof\ignorespaces}
\def\beginsubproof{
 \topsep6\p@\@plus6\p@\relax
\topsep4\p@\@plus4\p@\relax
\begin{trivlist}
  \item[\hskip\labelsep{
  \mbox{\ooalign
  {\hfill {\scriptsize $\because$}\hfill\crcr $\bigcirc$}}}]\rm}
\def\endsubproof{\end{trivlist}}
\makeatother                    
\makeatletter
\@addtoreset{claim}{section}
\@addtoreset{definition}{section}
\@addtoreset{theorem}{section}
\@addtoreset{lemma}{section}
\@addtoreset{remark}{section}
\@addtoreset{example}{section}
\@addtoreset{corollary}{section}
\@addtoreset{proposition}{section}
  

\begin{document}
\begin{center}
{\Large The Zero-Mass Limit Problem for a Relativistic Spinless Particle in an Electromagnetic Field 
} 
\end{center}

\begin{center} {\large
Takashi Ichinose   and Taro Murayama}
\end{center}
\begin{center}
Division of Mathematical and Physical Sciences, Graduate School of Natural Science and Technology, Kanazawa University Kakuma-machi, Kanazawa, 920-1192, Japan
\end{center}

\noindent{\bf Abstract:}   
It is shown that  mass-parameter-dependent
solutions of the imaginary-time magnetic relativitstic Schr\"odinger equations converge as functionals of  L\'evy processes represented by stochastic integrals of stationary Poisson point processes if mass-parameter goes to zero. 

\noindent {\bf 2010 Mathematics Subject Classification:}  60G51; 60F17; 60H05; 35S10; 81S40.

\noindent {\bf Key words:} magnetic relativistic Schr\"odinger operator, imaginary-time relativistic Schr\"odinger equation, L\'evy process, path integral formula, Feynman-Kac-It\^o formula. 
\section{Introduction and results.}
Kasahara-Watanabe \cite{K and W 86} discussed limit 
theorems in the framework of semimartingales represented by stochastic integrals of point processes. In fact, they considered  a sequence of point processes and their certain functionals represented by stochastic integrals, and proved  their convergence in that context.

In this paper we treat a sequence of a slightly more general functionals of special kind of L\'evy processes, which  have no  Gaussian part 
stemming from relativistic quantum mechanics,  to discuss its convergence.
Naturally we have in mind 
the following relativistic Schr\"odinger equation which describes a spinless quantum particle of mass $m>0$ (for example, pions)  in $\R^{d}$ under the influence of the vector and scalar potentials $A(x), V(x) \colon$
\begin{align}
&i\frac{\partial}{\partial t} \psi(x, t) = [H_{A}^{m} - m + V]\psi(x,t)
\quad  (t>0), 
\end{align}
where  $x\in\R^d$.
In this paper, to see the main idea, we only consider the case that
$A\in C^{\infty}_0 (\R^d;\R^d)$ and
$V  \in C_0(\R^d;\R)$. 
Here then $H_A^m$ is defined by
\begin{align*}
 &(H^{m}_{A}f)(x)  :=\text{Os-}\frac{1}{(2\pi)^{d} }   \iint_{\R^{d}\times \R^{d}} \nonumber 
e^{i(x-y)\cdot \xi} \sqrt{\left| \xi -A(\tfrac{x+y}{2})\right|^{2}+m^{2}}f(y)dy d\xi  \nonumber 
\end{align*}
for $f \in  C_{0}^{\infty}(\mathbf{R}^{d})$, where \lq \lq Os\rq \rq   means {\it oscillatory integral}. $H_{A}^{m}$ is called the {\itshape Weyl pseudo-differential operator with mid-point prescription}, corresponding to the classical relativistic Hamiltonian  $\sqrt{|\xi-A(x)|^{2}+m^{2}}$. 
It is  essentially selfadojoint in $L^{2}(\R^{d})$ on $C_{0}^{\infty}(\R^{d})$ and bounded from below by $m$ (\cite{I 89},\cite{I and T 93}). 
We have $H_{0}^{m}= \sqrt{-\Delta +m^{2}}$ for $A \equiv 0$, where $-\Delta$ is the Laplacian in $\R^{d}$.
The light velocity $c$, electric  charge $e$ and Planck's constant $h$ are taken to be $1$, $1$ and 2$\pi$ respectively.

The operator $H_{A}^{m}-m+V$ was first studied in \cite{I and T 86} by one of the authors of this paper to treat 
the {\it pure imaginary-time} relativistic  Schr\"odinger equation 
\begin{align}
&\frac{\partial}{\partial t} u(x, t) = -[H_{A}^{m} - m + V]u(x,t)
\quad   (t>0), 
\end{align}
where $x \in \R^d$. 
An imaginary-time path integral formula was given  on path space  $D_{0}$ to represent the solution of the Cauchy problem for (1.2).
Here $D_{0}$  
is the  set of the right-continuous paths $X: [0,\infty) \to  \mathbf{R}^{d}$ with left-hand limits and $X(0)=0$.

We use the probability space $(D_0,  \mathcal{F}, \lambda^{m})$  treated in  \cite{I and T 86} with the natural   filtration $\{\mathcal{F}(t)\}_{t \geq 0}$,  where $\mathcal{F}(t):=\sigma(X(s); s \leq t)\subset \mathcal{F}$. 
$\{X(t)\}_{t \ge 0}$ is L\'evy process, namely, it has stationary independent increments and is stochastically continuous (cf., \cite{I and W 81}, \cite{S 99}, \cite{A 09}). 
$\lambda^m(X; X(t) \in dy)$   is equal to   $k^{m}_{0}(y, t)dy$, where $k^{m}_{0}(y, t)$ is the integral  kernel of the operator $e^{-t(\sqrt{-\Delta +m^{2}}-m)}$ and 
has an explicit expression 
\begin{align}
     k_{0}^{m}(y,t) 
 =   \begin{cases}   
\displaystyle 2 \left( \frac{m}{2\pi}\right)^{(d+1)/2}\frac{te^{mt}  K_{(d+1)/2}(m(|y|^{2}+t^{2})^{1/2})}{(|y|^{2}+t^{2})^{(d+1)/4}},  \quad & m>0,
  \\  
 \displaystyle \frac{\Gamma ((d+1)/2)}{\pi^{(d+1)/2}} \frac{t}{(|y|^{2}+t^{2})^{(d+1)/2}}, & 
m=0.
 \end{cases}
\end{align}
Here $K_{\nu}$ stands for the modified Bessel function of the third kind of order $\nu$.
 
The characteristic function of $X(t)$   is
\begin{align}  
 & E^{m}[  e^{i  \xi \cdot X(t)}]=e^{-t(\sqrt{|\xi|^{2}+m^{2}}-m)}, \quad \xi \in \mathbf{R}^{d},
\end{align}
where  $E^{m}$  denotes the expectation over $D_0$ with respect to $\lambda^{m}$. 
By the {\itshape L\'evy-Khintchine formula},  
\begin{align}
 \sqrt{ |\xi|^{2} + m^{2}} - m  
 = - \int_{|y| > 0} \big( 
             e^{i\xi \cdot y} - 1 - i\xi  \cdot  y\mathbf{1}_{|y|<1} \big)
                n^{ m}   (dy). 
\end{align} 
Here $n^{m}(dy)$ is the {\itshape L\'evy measure}, that is  a $\sigma$-finite measure on $\mathbf{R}^{d}\setminus \{0\}$ satisfying $\int_{|y|>0}(1 \wedge |y|^{2})n^{m}(dy) < \infty$, and having  density 
\begin{align}\label{eq:eq levy measure}
n^{m}(y)=n^{m}(|y|)  
=  \begin{cases}
\displaystyle    2\left (\frac{m}{2\pi}\right)^{(d+1)/2}\frac{K_{(d+1)/2}(m|y|)}{|y|^{(d+1)/2}}, 
\quad & m>0, 
\\
 \displaystyle      \frac{\Gamma ((d+1)/2)}{\pi^{(d+1) /2}} \frac{1}{|y|^{d+1}}, & m=0.
\end{cases}
\end{align}
As shown in \cite{I 89}, $H_{A}^{m}$ has another expression  connected with the L\'evy measure  $n^{m}(dy)$
\begin{align*}
(H_{A}^{m}f)(x)=mf(x) -\lim_{r \downarrow 0}\int_{|y|\geq r} \big[e^{-iy\cdot A(x+\frac{1}{2}y)}  f(x+y) -f(x) \big]n^{m}(dy).
\end{align*}

For  $X \in D_0$,
let $N_{X}(dsdy)$ be  a counting measure on $(0, \infty) \times (\mathbf{R}^{d} \setminus \{0\})$ defined by 
\begin{align*}
&N_{X}(E):=\# \{ s >0; (s,X(s)-X(s-)) \in E\}
\end{align*}
for  $E \in \mathcal{B}(0, \infty) \times \mathcal{B}(\mathbf{R}^{d} \setminus \{0\})$, where $\mathcal{B}(\cdots)$ are   $\sigma$-algebras of Borel sets.
$N_{X}(dsdy)$ is the stationary Poisson random measure with   intensity measure $ds n^{m}(dy)$ with respect to $\lambda^{m}$. Let  $\widetilde{N_{X}^{m}}(dsdy):= N_{X}(dsdy)-ds  n^{ m}(dy)$. By the {\itshape L\'evy-It\^o theorem}, 
\begin{align}\label{eq:eq levy-ito}
X(t) =\int_{0}^{t}  \int_{|y| \geq 1}  y N_{X}(dsdy) 
+ \int_{0}^{t} \int_{0< |y| < 1}  y \widetilde{N_{X}^{m}}(dsdy),\quad \lambda^{m}\text{-a.s. }X \in D_{0}. 
\end{align}
Here and below, we should    understand   $\int_0^t: =\int_{(0, t]}$. 
It can be proved  that  the solution of (1.2)   with initial data $u^{m}(x,0)=g(x)$ is given by 
\begin{align}
u^{m}(x,t) &:= E^{m}[e^{-S^{m}(t,x,X)}g(x+X(t))],\\ 
 S^{m}( \cdot) & :=  i Y^m(t, x,X)  +  \int_{0}^{t}      V(x+X(s))ds,   \\
 Y^m (\cdot)  &:=  \int_{0}^{t}\int_{|y| \geq 1} A(x+X(s-)+\tfrac{1}{2}y) \cdot y N_{X}(dsdy) \nonumber \\
& \quad  + \int_{0}^{t}  \int_{0<|y| < 1}   A(x+X(s-)+\tfrac{1}{2}y) \cdot y \widetilde{N_{X}^{m}}(dsdy)\nonumber \\
& \quad  + \int_{0}^{t}  ds\int_{0<|y|<1}  \Big[ A(x+X(s)+\tfrac{1}{2}y)  -A(x+ X(s))\Big] \cdot y n^{m}(dy).     \nonumber   
\end{align}
In  (1.7) and (1.9) above,  the integration regions $|y| \geq 1$ and $0<|y| <1$ may be replaced by $|y| \geq \delta$ and $0< |y| < \delta$  respectively,  for any $\delta >0$.

We note that these relativistic quantities, $H^{m}_{A}-m+V$, $\sqrt{|\xi|^{2}+m^{2}}-m$, $D_{0}$, $\lambda^{m}$, $k_{0}^{m}(y,t)$ and  $X(t)$, correspond to the nonrelativistic ones $\frac{1}{2m}(-i \nabla -A)^{2}+V$, $\frac{|\xi |^{2}}{2m}$, $ C_{0}$,  {\itshape Wiener measure}, 
the heat kernel $(\tfrac{m}{2\pi t })^{d/2}e^{-\tfrac{m}{2t}|y|^{2}}$, {\itshape Brownian motion} $B(t)$, respectively. Here $C_{0}$ is the  space of  continuous paths $B: [0,\infty) \to  \mathbf{R}^{d}$ with $B(0)=0$. 
Furthermore, (1.8)  with (1.9)    is what {\it does}  correspond to {\itshape Feynman-Kac-It\^o formula} (\cite{S 79}).

The purpose of this paper is to answer the following question: \\
(Q) {\itshape When the mass $m > 0$ of the particle  becomes  sufficiently  small,  how does  its  property 
vary ?}

\begin{theorem} 
$\lambda^{m}$ converges weakly to $\lambda^{0}$ as $m \downarrow 0$.
\end{theorem}

\begin{theorem}
$u^{m}(\cdot, t)$ converges to  $u^{0}(\cdot,t )$ on  $L^{2}(\mathbf{R}^{d})$ as $m \downarrow 0$,  uniformly on  
$[0, T]$. 
\end{theorem}  

Here and below, $0 < T < \infty$ can be taken arbitrary.
Theorem 2 implies  the strong resolvent convergence of $H^{m}_{A}-m+V$ to $H^{0}_{A}+V$ ([13, IX, Theorem 2.16]). An immediate consequence is the following result   for the solution   $\psi^{m}(x,t)$ of the Cauchy problem for (1.1).

\noindent
\begin{corollary}
$\psi^{m}(\cdot,t)$ converges to  $\psi^{0}(\cdot ,t)$ on  $L^{2}(\mathbf{R}^{d})$ as $m \downarrow 0$, uniformly on $[0, T]$.
\end{corollary}

We will prove  Theorem 2 by using following: 

\begin{theorem}
 $u^{m}(\cdot, t)$ converges to  $u^{0}(\cdot,t )$ on  $C_{\infty}(\mathbf{R}^{d})$ as $m \downarrow 0$, uniformly on $[0, T]$, where $C_{\infty}(\mathbf{R}^{d})$  is the space of the  continuous functions $g :\mathbf{R}^{d} \to \mathbf{C}$ with $|g(x)| \to 0$ as $|x| \to \infty$ with norm $\|g\|_{\infty}:= \sup_{x\in \R^{d} }|g(x)|$.
\end{theorem} 

The crucial idea of proof is to do a change of variable \lq \lq path\rq\rq.
In Sections 2,3 and 4,  these theorems are shown  by probabilistic method, although one can more easily show Theorem 2 by operator-theoretical one \cite{I 92},  and also by pseudo-differential calculus \cite{N and U 90}.  
In this paper, as we mentiond before, we treat the problem under 
a rather mild assumption on the  potentials $A(x), V(x)$. 
We will come to more general case in a forthcoming paper, together for  the other two different magnetic relativistic Schr\"odinger operators  (\cite{I 2012}, \cite{I 2013}) corresponding to the same classical relativistic  Hamiltonian. Another limit problem when  the light velocity $c$ goes to infinity ({\it nonrelativistic limit}) was studied in \cite{I 87}.

\section{Proof of Theorem 1.}  
We observe the following three facts which imply Theorem 1 ([2,Theorem 13.5]): \\
(i) The finite dimensional distributions with respect to $\lambda^{m}$ converge  weakly to those with respect to $\lambda^{0}$ as $m \downarrow 0$.\\ 
(ii) For any $t > 0$,   $\lambda^0(X; X(t)-X(t-\varepsilon)\in dy )$ converges weakly to  Dirac measure concentrated at the point $0\in \R^d$ as $\varepsilon    \downarrow 0$. \\ 
(iii) There exist  constants $\alpha > \frac{1}{2}$, $\beta > 0$ and a nondecreasing continuous function $F$ on $[0,\infty)$  such that 
\begin{align*}
& E^{m}\left[ |X(s)-X(r)|^{\beta} |X(t)-X(s)|^{\beta}\right]  \leq  [F(t)-F(r)]^{2\alpha}, \quad 0 < m < 1, \ 0 \leq r < s < t.
\end{align*}
 {\itshape Proof.}  (i) follows from (1.4), and (ii) from the stochastic continuity of $\{X(t)\}_{t \geq 0}. $     
(iii) Since $\tfrac{d}{d\tau} \tau^{\nu} K_{\nu}(\tau)=$ $-\tau^{\nu}K_{\nu-1}(\tau) $ ($\tau>0, \nu >0$) ([3,  (21), p.79]) and $\nu \mapsto K_{\nu}(\tau)$ is strictly increasing in $(0, \infty)$  ([3, (21), p.82]), we have $(d/d\tau)(e^{\tau}\tau^{\nu}K_{\nu}(\tau))= e^{\tau}\tau^{\nu} (K_{\nu}(\tau)-K_{\nu-1}(\tau))<0$  if $0< \nu <\frac{1}{2}$.  Therefore    $\tau \mapsto e^{\tau}\tau^{\nu} K_{\nu}(\tau)$ is strictly decreasing in $(0,\infty)$ and so [3, (41),(42), (43), p.10]
\begin{align}
e^{\tau} \tau^{\nu}K_{\nu}(\tau)   \leq  \lim_{\tau \downarrow 0} \tau^{\nu}K_{\nu}(\tau)
                                   =2^{\nu-1} \Gamma(\nu).
\end{align}
Then we have for $0\leq r<  s< t$, $\frac{1}{2}< \beta <  1$, 
\begin{align*}
 E^{m}\left[ |X(s)-X(r)|^{\beta} |X(t)-X(s)|^{\beta}\right]
& = \int|y|^{\beta}k_{0}^{m}(y, s-r)dy   \int|y|^{\beta}k_{0}^{m}(y, t-s)dy\\
&= C(d, \beta)^{2} ((s-r)(t-s) )^{\beta} \\
&  \quad  \times e^{m(s-r)}   (m(s-r))^{ \frac{1-\beta}{2 }}K_{\frac{1-\beta}{2}} (m(s-r))\\
&  \quad  \times e^{m(t-s)}   (m(t-s))^{ \frac{1-\beta}{2 }}K_{\frac{1-\beta}{2}} (m(t-s))\\
&\leq  C(d, \beta)^{2}    2^{-(1+2\beta)}\Gamma(\tfrac{1-\beta}{2})^{2}(t-r)^{2\beta},
\end{align*}
where in the second equality we use [4, Lemma 3.3(ii)] with a constant  $C(d, \beta)$ depending on $d$ and $\beta$.
Therefore  (iii) holds for $\frac{1}{2}< \beta < 1$ and  $\alpha = \beta$ and $F(p):= C(d, \beta)^{1/\beta} 2^{-(1+2\beta)/2\beta}\Gamma(\tfrac{1-\beta}{2})^{1/\beta}p$.\hfill$\square$

\section{Proof of Theorem 2.} We will prove Theorem 2  by  assuming   validity of Theorem 3. In this and the next section, we assume  $V \geq 0$ without loss of generality,  since  in the general case, we have only to replace $V$ in (1.8), (1.9) by $V -\inf V \geq 0$.\\ 
Step I:  Let $g \in C_{0}^{\infty}(\R^{d})$.  For $R >0$, we have
\begin{align*}
\|u^{m}(\cdot, t )-u^{0}(\cdot, t )\|_{2}
&\leq  \|u^{m}(\cdot, t )-u^{0}(\cdot, t )\|_{L^{2}(|x|< R)}  +\|u^{m}(\cdot, t )-u^{0}
(\cdot, t )\|_{L^{2}(|x|\geq R)}\nonumber \\
{} &=: I_{1}( t, m,R)+I_{2}( t, m,R).
\end{align*}
From  Theorem 3, $I_{1}( t, m,R)$ converges to zero as $m \downarrow 0$ uniformly on $t \leq T$.
From (1.8), we have 
\begin{align*}
I_{2}(t, m,R)& \leq  \|u^{m}(\cdot,t)\|_{L^{2}(|x| \geq R)} +      \|   u^{0}(\cdot,t)\|_{L^{2}(|x| \geq R)} \nonumber \\
&\leq   \Big(\int_{|x| \geq R}dx \int k_{0}^{m}(y,t)|g(x+y)|^{2}dy\Big)^{\frac{1}{2}}  \\
& \quad  + \Big(\int_{|x| \geq R} dx \int k_{0}^{0}(y,t)|g(x+y)|^{2}dy\Big)^{\frac{1}{2}}  \\
&=:  J( t, m,R)+J(t,0,R).
\end{align*}
 Let $\chi$ be a nonnegative $C^{\infty}_{0}(\mathbf{R}^{d})$ function such that $\chi(x)=1$ if $|x| \leq \frac{1}{2}$ and $=0$ if $|x| \geq 1$. Put $h(x)=|g(x)|^{2}$. Since $\mathbf{1}_{|x|< R} \geq \chi(\frac{x}{R})$, we have
\begin{align*}
J( t, m,R)^{2} 
& \leq  \int(1-\chi(\tfrac{x}{R})  )dx\int k_{0}^{m}(y,t)h(x+y)dy  \\
& = \frac{1}{(2\pi)^{d}} \bigg [\widehat{h}(0)  \int \left(  1- \exp \{ -t[\sqrt{\tfrac{|\eta|^{2}}{R^{2}}+m^{2}}-m] \}  \right)  \overline{\widehat{\chi}(\eta)}  d\eta  \\
&  \qquad \qquad  +    \int (  \widehat{h}(0)) - \widehat{h}(\tfrac{\eta}{R})) \exp \left\{ -t\left[\sqrt{\tfrac{|\eta|^{2}}{R^{2}}+m^{2}}-m\right] \right\}\overline{\widehat{\chi}(\eta)} d\eta  \bigg],
\end{align*}
which converges to zero as $R \to \infty$ uniformly on $t \leq T$ and $0\leq m  \leq1$. Here,  for $\varphi \in \mathcal{S}(\R^{d})$, $\widehat{\varphi}$ is the Fourier transform of $\varphi$ given by $\widehat{\varphi}(\xi)=\int e^{-ix \cdot \xi}\varphi(x)dx$ ($\xi \in \R^{d}$).

From  (1.3) and (2.1), it follows that $k^{m}_{0}(y,t) \to k_{0}^{0}(y,t)$ as $m \downarrow 0$, and then 
$ J( t, 0,R)^{2}\leq \liminf_{m \downarrow 0} J( t, m,R)^{2}$ by Fatou's lemma. Therefore we have Theorem 2 for this step. \\ 
Step II:  Let  $g \in L^{2}(\mathbf{R}^{d})$.  There is a sequence $\{ g_{n}\} \subset C^{\infty}_{0}(\mathbf{R}^{d})$  such that 
$g_{n} \to g$ in $L^{2}(\mathbf{R}^{d})$ as $n \to \infty.$ Put $u^{m}_{n}(x,t):=  E^{m}[e^{-S^{m}(t,x,X)}g_{n}(x+X(t))]$.     Then we have
\begin{align*}
\|   u^{m}(\cdot, t )-u^{0}(\cdot, t )\|_{2} 
&\leq     \|   u^{m}(\cdot, t )-u^{m}_{n}(\cdot, t )\|_{2}+  \|   u^{m}_{n}(\cdot, t )-u^{0}_{n}(\cdot, t )\|_{2}\\ 
& \quad +  \|   u^{0}_{n}(\cdot, t )-u^{0}(\cdot, t )\|_{2}\\ 
&\leq 2\|g_{n}-g\|_{2}+\| u^{m}_{n}(\cdot, t )-u^{0}_{n}(\cdot, t ) \|_{2}.
\end{align*}
By Step I, we have
\begin{align*}
\limsup_{m \downarrow 0}\sup_{t \leq T}\|    u^{m}(\cdot, t )-u^{0}(\cdot, t )                  \|_{2} \leq 2 \|g_{n}-g\|_{2},
\end{align*} which converges  to zero as $n \to \infty$.  \hfill$\square$

\section{Proof of Theorem 3.} 
From  (1.8), we have to  prove that 
\begin{align*}
u^{m}(x,t)  &= E^{m}[e^{-S^{m}(t,x,X)}g(x+X(t))] \nonumber   \\
                             & \to E^{0}[e^{-S^{0}(t,x,X)}g(x+X(t))] = u^{0}(x,t)
\end{align*}
as  $m \downarrow 0$ in $C_{\infty}(\mathbf{R}^{d})$. But its direct proof seems difficult since both  the integrand $e^{-S^{m}(t,x,X)}g(x+X(t))$ and the probability measure $\lambda^{m}$  depend on $m$.  So we change   $E^{m}[\cdots]$   
to  $E^{0}[\cdots]$  by a  {\itshape change of variable} (i.e., {\it change of probability measure}) $\lambda^{m}=\lambda^{0} \Phi_{m}^{-1}$ with path space transformation 
$\Phi_{m}: D_0 \to D_0 $. 
If there is  such a $\Phi_{m}$,  we can see by (1.4) and (1.5) that  
the difference between the path $X(t)$ and the transformed path $\Phi_{m}(X)(t)$   is expressed  
in terms of  the difference between the two L\'evy measures $n^{0}(dy)$ and $n^{m}(dy)$, so that   it is  presumed to hold that $n^{m}(dy)=n^{0}\phi_{m}^{-1}(dy)$   
for some map $\phi_{m}\colon {\mathbf{R}^{d}\setminus \{0\} }$ $\to$ $\mathbf{R}^{d}\setminus \{0\}$.

We will determine  $\phi_{m}$  in such a way that (1) $n^{m}(dy)=n^{0}\phi_{m}^{-1}(dy)$, (2)  $\phi_{m} \in C^{1}(\mathbf{R}^{d}\setminus \{0\}; \mathbf{R}^{d}\setminus \{0\})$, (3) $\phi_{m}$ is one to one and onto, (4) $\det D\phi_{m}(z) \ne 0$ for all $z \in \R^{d}\setminus \{0\}$,  where $D\phi_{m}(z)$ is the Jacobian matrix
of $\phi_{m}$ at the point $z$. 

Let $U:= \{y \in \R^{d} \setminus \{0\}; |y| \in U'\}$ for $U' \in \mathcal{B}(0, \infty)$.
Introducing the spherical coordinates by $z=r \omega$, $r>0$, 
$\omega \in S^{d-1}$, we have
\begin{align*}
n^{m}(U)=\int_{U}n^{m}(|y|)dy 
= C(d)\int_{U'}n^{m}(r)r^{d-1}dr,
\end{align*}
where $C(d)$ is the surface area of the $d$-dimensional unit ball. 

Let us assume that $\phi_{m}^{-1}(z) = l_{m}(|z|)\tfrac{z}{|z|}$ for some non-decreasing $C^{1}$ function $l_{m}:(0, \infty) \to (0, \infty)$. 
Then we have
\begin{align*} 
n^{0}\phi_{m}^{-1}(U) 
&= \int_{U}n^{0}(l_{m}(|z|))   |z|^{-(d-1)} l_{m}(|z|)^{d-1}l_{m}'(|z|)     dz\nonumber \\
&  = C(d)\int_{U'} n^{0}(l_{m}(r))   l_{m}(r)^{d-1} l_{m}'(r)  dr,
\end{align*}
where $l_{m}'(r)=(d/dr)l_{m}(r).$ Therefore we have 
\begin{align*}
n^{m}(r)r^{d-1}= n^{0}(l_{m}(r))   l_{m}(r)^{d-1} l_{m}'(r),       \quad \text{a.s. } r > 0.  
\end{align*}
If $m > 0$, from (1.6),   we have
\begin{align*}
&-\frac{d}{dr}l_{m}(r)^{-1}= 2^{-\frac{d-1}{2}} \Gamma(\tfrac{d+1}{2})^{-1}  m^{\frac{d+1}{2}}r^{\frac{d-3}{2}} K_{\frac{d+1}{2}}(mr).
\end{align*}
We solve this differential equation under  boundary condition $l_{m}(\infty) = \infty$ to    get 
\begin{align}
l_{m}(r)=  \frac{2^{\frac{d-1}{2}} \Gamma (\frac{d+1}{2})}{ m^{\frac{d+1}{2}}   \int_{r}^{\infty}    u^{\frac{d-3}{2}} K_{\frac{d+1}{2}}(mu)du }.
\end{align}
Here we note that $0<\int_{r}^{\infty}   u^{\frac{d-3}{2}} K_{\frac{d+1}{2}}(mu) du < \infty$  by  $K_{\frac{d+1}{2}}(\tau) > 0$ for $\tau > 0$, and [3, (37), (38), p.9]
\begin{align*}
K_{\frac{d+1}{2}}(\tau)    =\left(\frac{\pi}{2}\right)^{1/2} \tau^{-1/2}e^{-\tau}  (1+o(1)), \quad  \tau \uparrow \infty.
\end{align*}
\begin{proposition}
{\rm (i)} $l_{m}(r)$ is a strictly increasing  $C^{\infty}$ function of $r\in(0, \infty)$ and $l_{m}(+0)=0$, $l_{m}(\infty)=\infty$.\\
 {\rm (ii)} For all $r>0$, $l_{m}(r)$ converges to $r$, strictly decreasingly,  as  
$m \downarrow 0$.  
\end{proposition}

\noindent
{\itshape Proof.}  (2.1) implies $l_{m}(+0)=0$. The other claims of (i) follow from (4.1) and the fact that  $K_{(d+1)/2}(\tau)$ is a $C^{\infty}$ function in $(0,\infty)$. The claim (ii) can be proved by  the fact that $\tau^{\nu}K_{\nu}(\tau)$ is strictly decreasing in $(0, \infty)$ (cf. Section 2, Proof of (ii)), (2.1) and the monotone convergence theorem.   \hfill$\square$

If $m=0$, let  $l_{0}(r):=r$. Let  us put    $\phi_{0}(z):=z$ and for $m >0$,  
\begin{align*}
\phi_{m}(z):= \l_{m}^{-1}(|z|) \frac{z}{|z|}, \quad z \in \R^{d}\setminus \{0\}.
\end{align*}
Then we have
\begin{align*}
\phi_{m}^{-1}(z)= \l_{m}(|z|) \frac{z}{|z|}, \quad z \in \R^{d}\setminus \{0\}.
\end{align*}
We note that  
\begin{align}
\phi_{m}(z) \to z, \; |\phi_{m}(z)|= l_{m}^{-1}(|z|) \uparrow |z|
\end{align} \text{as} $m \downarrow 0$  by Proposition 1 (ii).

Let us define  $\Phi_{0}(X):=X$ and for $m >0$, 
\begin{align}
\Phi_{m}(X)(t) 
& :=  \int_{0}^{t}\int_{|y| \geq 1}y N_{X}(ds \phi_{m}^{-1}(dy))    + \int_{0}^{t} \int_{0< |y| < 1} y \widetilde{N_{X}^{0}}(ds \phi_{m}^{-1}(dy))  \nonumber \\
&  = \int_{0}^{t} \int_{|z| \geq l_{m}(1)} \phi_{m}(z) N_{X}(ds dz)+\int_{0}^{t}\int_{0< |z| < l_{m}(1)}   \phi_{m}(z) \widetilde{N_{X}^{0}}( ds  dz)\nonumber \\
&  = \int_{0}^{t}  \int_{|z| \geq 1}\phi_{m}(z) N_{X}(ds dz)+ \int_{0}^{t}\int_{0< |z| < 1}   \phi_{m}(z) \widetilde{N_{X}^{0}}( ds  dz).
\end{align}

\begin{proposition} For every sequence $\{m\}$ with $m \downarrow 0$, there exists a subsequence $\{m'\}$ such that 
\begin{align*}
\sup_{t \leq T}|\Phi_{m'}(X)(t)-X(t)|\to 0  \text{ as } m' \downarrow 0, \lambda^{0}\text{-a.s. } X \in    D_0.
\end{align*}
\end{proposition}
\noindent
{\itshape Proof.} 
From (1.7) and (4.3), we have
\begin{align*}
 \sup_{t \leq T}|\Phi_{m}(X)(t)-X(t)|\nonumber   & \leq \int_{0}^{T}\int_{|z| \geq 1}| \phi_{m}(z)-z|N_{X}(dsdz)\\ 
&  \quad +  \sup_{t \leq T}\Big| \int_{0}^{t}\int_{0< |z| <1} (\phi_{m}(z)-z) \widetilde{N_{X}^{0}}(dsdz)\Big|   \nonumber\\
&=: I_{1}(m, X)+ \sup_{t \leq T}|I_{2}(t, m, X)|.
\end{align*} 
We have  $I_{1}(m, X) \to 0$ as $m \downarrow 0$ by  (4.2) and 
$\int_{0}^{T}\int_{|z|\geq 1}|z|N_{X}(dsdz) 
< \infty$.
We note  that $I_2(t,m,X)$ is  the 
$L^2(D_0;\lambda^0)$-limit of the 
right-continuous $\{\mathcal{F}(t)\}_{t \geq 0}$-martingale 
 $\{I_2^{\varepsilon}(t,m,X)\}_{t \geq 0}$ with  
$I_2^{\varepsilon}(t,m,X):= \int_0^t\int_{\varepsilon<|z|<1}(\phi_m(z)-z)\widetilde{N_X^0}(dsdz)$  
as $\varepsilon \downarrow 0$, with convergence being uniform on $t \leq T$.  By taking a subsequence if necessary,   $I_2^{\varepsilon}(t,m,X)$ converges to $I_2(t,m,X)$ as $\varepsilon  \downarrow  0$ uniformly on $t \leq T$, 
 $\lambda^0$-a.s., and hence $I_2(t,m,X)$ is right-continuous on $t \leq T$,    $\lambda^0$-a.s.        
([11, p.73, Proof of Theorem 5.1], [15, p.128-129, Proofs of Lemmas 20.6, 20.7]).
Then we use 
Doob's martingale inequality \cite{A 09} to have 
\begin{align*}
E^{0}\left[ \sup_{t \leq T}|I_{2}(t, m, X)|^{2} \right]
 &\leq  4E^{0}\left[| I_{2}(T, m, X)|^{2} \right]\nonumber    \\
& \leq 4T \! \! \int_{0<|z|<1}\!\! \!\!\!\!\!\!\!\!\! \!|   \phi_{m}(z)-z|^{2} n^{0}(dz),\nonumber 
\end{align*}
which converges to zero as $m \downarrow 0$ by (4.2) and $\int_{0< |z|<1} |z|^{2}n^{0}(dz) < \infty$. \hfill$\square$

By  (1.8) and $\lambda^{m}=\lambda^{0} \Phi_{m}^{-1}$, we have 
\begin{align*}
u^{m}(x,t)=   E^{0}[  e^{-S^{m}(t, x, \Phi_{m}(X))} g(x+ \Phi_{m}(X)(t))],
\end{align*} 
and then 
\begin{align}
\sup_{t \leq T}   \|u^{m}(\cdot,t) -u^{0}(\cdot,t)\|_{\infty}
 \leq  &  \| g \|_{\infty} \sup_{t \leq T, \ x \in \mathbf{R}^{d}} E^{0} \left[\left| e^{-S^{m}(t,x, \Phi_{m}(X))}  -  e^{-S^{0}(t, x, X)} \right|   \right] \nonumber \\ 
&  +E^{0}\left[ \sup_{t\leq T} \|g(\cdot + \Phi_{m}(X)(t))  -  g(\cdot +  X(t))\|_{\infty}  \right].\end{align}      
Since $g\in C_{\infty}(\R^{d})$ is uniformly continuous and bounded on $\R^{d}$,  the second term on the right  of (4.4)  converges to zero  as $m \downarrow 0$. 

Next we consider the first term on the right  of (4.4). 
By $N_{\Phi_{m}(X)}(dsdy) = N_{X}(ds \phi_{m}^{-1}(dy))$,
we have 
\begin{align}\label{eq;eq S}
&S^{m}(t, x, \Phi_{m}(X)) \nonumber \\
&=  i\bigg(\int_{0}^{t}  \int_{|z|\geq 1} A(x+ \Phi_{m}(X)(s-)+\tfrac{1}{2}\phi_{m}(z))\cdot \phi_{m}(z) N_{X}(dsdz)  \nonumber \\
 &\quad + \int_{0}^{t}  \int_{0< |z|<1} A(x+ \Phi_{m}(X)(s-)+\tfrac{1}{2}\phi_{m}(z))\cdot \phi_{m}(z) \widetilde{N_{X}^{0}}(dsdz)  \nonumber \\
& \quad +\int_{0}^{t}  ds  
  \int_{0<|z| < 1}  \left[A(x+\Phi_{m}(X)(s)+\tfrac{1}{2}\phi_{m}(z))-A(x+\Phi_{m}(X)(s))\right]  \cdot \phi_{m}(z) n^{0}(dz) \bigg)  \nonumber   \\
&\quad +    \int_{0}^{t}  V(x+ \Phi_{m}(X)(s))ds   \nonumber \\
& =:  i \Big(  S_{1}^{m} (t,x,X)+S_{2}^{m}(t,x,X)+ S_{3}^{m} (t,x,X) \Big)     + S_{4}^m(t,x,X).  \nonumber 
\end{align}
 By the inequality
\begin{equation*}
 |e^{-(ia+b)}-e^{-(ia'+b')}| \leq e^{-b}|e^{-ia}-e^{-ia'}|+|b-b'|
\end{equation*}  for any $a,a'\in\mathbf{R}$, $b, b'  \geq 0$, 
$\sup E^{0}[\cdots]$ of  the first term on the right of  (4.4) is less than or equal to 
\begin{align}  
 &   E^{0}\Big[\sup_{t \leq T} \|e^{-i S^m_{1}(t, \cdot, X)}- e^{-i S^0_{1}(t, \cdot,X)}  \|_{\infty} \Big]    +  \sup_{x \in \mathbf{R}^{d}}E^{0}\Big[ \sup_{t \leq T} | S^m_{2}(t, x,  X)- S^0_{2}(t, x,X)  |\Big] \nonumber \\
 &     \quad  +E^{0}\Big[  \sup_{t \leq T} \| S^m_{3}(t, \cdot, X)-S^0_{3}(t, \cdot,X)    \|_{\infty}  \Big]  + E^{0}\Big[\sup_{t \leq T}  \| S^m_{4}(t, \cdot, X) -S^0_{4}(t, \cdot,X)  \|_{\infty} \Big]. 
\end{align}
Now, let $\{m\}$ be a sequence with $m \downarrow 0$ and $\{m'\}$  any subsequence of $\{m\}$.  By Proposition 2, there exists a subsequence $\{m''\}$ of $\{ m'\}$ such that  $\sup_{t \leq T}|\Phi_{m''}(X)(t)-X(t)|\to 0$  \text{ as } $m'' \downarrow 0$, $\lambda^{0}$\text{-a.s}. 

To prove that each term of (4.5) converges  to zero as $m'' \downarrow  0$, we first  note that 
\begin{align*}
   S^{m''}_{1}(t, x, X) -S_{1}^0(t, x,X)  & = \int_{0}^{t} \int_{|z| \geq  1 }
\big( A(x+\Phi_{m''}(X)(s-) + \tfrac{1}{2}\phi_{m''}(z))  \nonumber \\
&\qquad \qquad \qquad  -  A(x+X(s-)+\tfrac{1}{2}z)\big)\cdot  \phi_{m''} (z) N_{X}(dsdz)\\
&  \quad +\int_{0}^{t} \int_{|z| \geq  1 } A(x+X(s-)+\tfrac{1}{2}z)\cdot   (\phi_{m''}(z)- z)  N_{X}(dsdz).
\end{align*}
Then the integrand of the first term of (4.5) is less than or equal to  
\begin{align*} 
&  \!\int_{0}^{T} \int_{|z| \geq  1}    \sup_{x\in \R^{d}}\left|A(x\!+\!\Phi_{m''}(X)(s-) + \tfrac{1}{2}\phi_{m''}(z))   -  A(x+X(s-)+\tfrac{1}{2}z)\right|  |z| N_{X}(dsdz)\\
&  +\sup_{x \in \R^{d}}|A(x)|\!\int_{0}^{T}\!\! \int_{|z| \geq  1 }\!\!\!\!\!   |\phi_{m''}(z)\!-\! z|  N_{X}(dsdz),
\end{align*}
which converges to zero as $m'' \downarrow 0$ since $A \in C_{0}^{\infty}(\R^{d};\R^{d}) $  is uniformly continuous on $\R^{d}$.

Next, since $S_2^m(t,x,X)$ is seen  to be right-continuous, 
 by Schwarz's inequality and Doob's martingale inequality, $E^{0}[\cdots]$ of 
the second term of (4.5) is less than or equal to
\begin{align*}
&  2E^{0}\Big[ \int_{0}^{T} ds \int_{0<|z| <1}  \big|A(x\!+\!\Phi_{m''}(X)(s-)+\tfrac{1}{2} \phi_{m''}(z))\cdot\phi_{m''}(z)\\
 &\qquad \qquad \qquad  \qquad \quad -  A(x+X(s-)+\tfrac{1}{2}z) \cdot z \big|^{2}           n^{0}(dz)  \Big  ]^{\frac{1}{2}}.\nonumber 
\end{align*}   
By  the inequality  $( a+ b)^{2}\leq 2(a^{2}+b^{2})$ for any  $a,b \in \R$,   $E^{0}[\cdots]$ above  is less than or equal to
\begin{align*}
& 2\bigg\{E^{0}\bigg[ \int_{0}^{T}ds \int_{0<|z| <1} \sup_{x\in \R^{d}} |A(x+\Phi_{m''}(X)(s-) +\tfrac{1}{2} \phi_{m''}(z)) \\
& \qquad \qquad \qquad \qquad \qquad \qquad   -   A(x+X(s-)+\tfrac{1}{2} z)  |^{2} |z|^{2}  n^{0}(dz)\bigg]    \nonumber \\
& \quad +  T \sup_{x\in \R^{d}}|A(x)|^{2}\int_{0<|z| <1}    |\phi_{m''}(z)-z|^{2} n^{0}(dz)\bigg\},  \nonumber 
\end{align*}
which converges to zero as $ m'' \downarrow 0$.
As for the third term of (4.5), by the mean  value theorem, we have
\begin{align*}
& S^{m''}_{3}(t, x, X)-S^{0}_{3}(t, x,X)  \\
&=  \frac{1}{2}\int_{0}^{t}ds   \int_{0<|z| <1}     n^{0}(dz)\int_{0}^{1} \big[ \big(W^{m''}_{x,X}(s,\theta) \phi_{m''}(z) \big)\cdot  \phi_{m''}(z)   - \big( W^{0}_{x,X}(s,\theta)  z \big) \cdot  z \big]
d\theta.\end{align*}
Here $W^{m''}_{x,X}(s,\theta) $ and $W^{0}_{x,X}(s,\theta) $ are 
$d\times d$ matrices defined by   
\begin{align*}
&W^{m''}_{x,X}(s,\theta) = DA (x+ \Phi_{m''}(X)(s)+\tfrac{1}{2}\phi_{m''}(z) \theta),\\
&W^{0}_{x,X}(s,\theta) = DA (x+ X(s)+\tfrac{1}{2}z \theta),
\end{align*}
where $DA(\cdot )$ is the Jacobian matrix of $A$. 
Since 
\begin{align*}
&\big(W^{m''}_{x,X}(s,\theta)  \phi_{m''}(z) \big) \cdot  \phi_{m''}(z)
-  \big(W^{0}_{x,X}(s,\theta) z \big) \cdot z \\
&= \big(W^{m''}_{x,X}(s,\theta) \phi_{m''}(z) \big) \cdot 
     \big(\phi_{m''}(z) - z  \big) + \big( (W^{m''}_{x,X}(s,\theta)-W^{0}_{x,X}(s,\theta))\phi_{m''}(z) \big) \cdot  z \\
&\quad  + \big( W^{0}_{x,X}(s,\theta)(\phi_{m''}(z) - z) \big) \cdot z , 
\end{align*}
the integrand of the third term of (4.5) is less than or equal to
\begin{align*}
&T  \sup_{x \in \R^{d}}\| DA(x)\| \int_{0< |z| <1}   |\phi_{m''}(z) -z|  |z| n^{0}(dz)\\& +\frac{1}{2} \int_{0}^{T}ds \int_{0 < |z| <1} |z|^{2} n^{0}(dz)  \int_{0}^{1} \sup_{x \in \R^{d}} \| W^{m''}_{x,X}(s,\theta)-  W^{0}_{x,X}(s,\theta)\|d\theta,\nonumber 
\end{align*}
where $\| \cdot \|$ is the norm of matrices. 
This is less than or equal to 
\begin{align*}
3T  \sup_{x\in \R^{d}}\| DA(x)\| \int_{0< |z| <1}   |z|^{2} n^{0}(dz) < \infty,
\end{align*}
and 
converges to zero as $m'' \downarrow 0$ 
because each component of $DA$ is uniformly continuous on $\R^{d}$.

Finally, 
the fourth term of (4.5) is less than or equal to
\begin{align*}
E^{0}\left[\int_{0}^{T}\|V(\cdot+ \Phi_{m''}(X)(s))- V(\cdot+ X(s))\|_{\infty}ds\right],\end{align*}
 which converges  to zero as $m'' \downarrow 0$ since $V \in C_{0}(\R^{d};\R )$ is uniformly continuous on $\R^{d}$. 
Thus  we have  
$\sup_{t \leq T}   \|u^{m''}(\cdot,t) -u^{0}(\cdot,t)\|_{\infty}\to 0$ as $m'' \downarrow 0$,
and hence 
$\sup_{t \leq T}   \|u^{m}(\cdot,t) -u^{0}(\cdot,t)\|_{\infty}\to 0$ as $m \downarrow 0$.
\hfill$\square$

\vspace{2mm}           
\noindent
\textbf{Acknowledgements.} 
The  author (T. I.) is grateful to Professor Yuji Kasahara for a number of helpful discussions with suggestion on the subject at a very early stage of this work and Professor Masaaki Tsuchiya for frequent valuable and helpful  discussions from the beginning.
The other author (T.M.) would like to thank Professor Hidekazu  
Ito for  his   kind guidance,  many helpful advices and 
 warm  encouragement during the preparation of this work. 
The authors are indebted to the anonymous referee for valuable comments 
 and suggestions.


\begin{thebibliography}{16}
\bibitem{A 09} D. Applebaum, {\itshape L\'evy Processes and Stochastic Calculus}, 2nd ed.,  Cambridge Studies in Advanced Mathematics, 116, Cambridge Univ. Press, Cambridge, 2009.
\bibitem{B 99} P. Billingsley,  {\itshape Convergence of Probability Measures}, 2nd ed., Wiley Series in Probability and Statistics: Probability Statistics,  New York, 1999.
\bibitem{E 53} A. Erd\'elyi, W. Magnus, F. Oberhettinger and F. G. Tricomi, {\itshape  Higher Transcendental Fnctions}, {\it Vol. II}, McGraw-Hill, New York, 1953.
\bibitem{I 87}T. Ichinose,  The nonrelativistic limit problem for a relativistic spinless particle in an electromagnetic field,  J. Funct. Anal. \textbf{73} (1987), no. 2,  233--257.
\bibitem{I 89}T. Ichinose,  Essential selfadjointness of the Weyl quantized relativistic Hamiltonian,  Ann. Inst. H. Poincar\'e, Phys. Th\'eor. {\bf 51} (1989), no.3, 265--297. 
\bibitem{I 92} T. Ichinose, Remarks on the Weyl quantized relativistic Hamiltonian, Note Mat.  \textbf{12} (1992), 49--67.
\bibitem{I 2012}T. Ichinose, On three magnetic relativistic Schr\"odinger  operators and imaginary-time path integrals, Lett. Math. Phys \textbf{101} (2012), 323--339.
\bibitem{I 2013}T. Ichinose,  Magnetic relativistic Schr\"odinger operators and imaginary-time path integrals, in  {\it Mathematical Physics, Spectral Theory and Stochastic Analysis}, Oper. Theory Adv. Appl.,  232,  Birkh\"auser/Springer Basel AG, Basel,  2013, pp.247--297.   
\bibitem{I and T 86} T. Ichinose and H. Tamura, Imaginary-time path integral for a relativistic spinless particle in an electromagnetic field,  Commun. Math. Phys. \textbf{105} (1986), no. 2,  239--257.  
\bibitem{I and T 93} T. Ichinose and T. Tsuchida,  On essential selfadjointness of the Weyl quantized relativistic Hamiltonian,  Forum Math. {\textbf 5} (1993), no. 6, 539--559. 
\bibitem{I and W 81} N. Ikeda and S. Watanabe, {\itshape Stochastic Differential Equations and Diffusion Processes}, North-Holland Mathematical Library, 24, North-Holland, Amsterdam, 1981.
\bibitem{K and W 86} Y. Kasahara and S. Watanabe,  Limit theorems for point processes and their functionals,  J. Math. Soc. Japan {\textbf  3\textbf{8}} (1986), no. 3, 543--574.
\bibitem{K 76} T. Kato, {\itshape Perturbation Theory for Linear Operators}, 2nd ed., Springer, Berlin, 1976. 
\bibitem{N and U 90} M. Nagase and T. Umeda,  Weyl quantized Hamiltonians of relativistic spinless  particles in magnetic fields,  J. Funct. Anal. \textbf{92} (1990), no. 1, 136--154. 
\bibitem{S 99} K. Sato, {\itshape L\'evy Processes and Infinitely Divisible Distributions}, translated from the 1990 Japanese original,  Cambridge Studies in Advanced Mathematics, 68, Cambridge Univ. Press, Cambridge, 1999. 
\bibitem{S 79} B. Simon, {\itshape Functional Integration and Quantum Physics}, Pure and Applied Mathematics, 86, Academic Press, New York, 1979; 2nd ed.,  AMS Chelsea Publishing, Providence, RI, 2005.
\end{thebibliography}
\end{document}